\documentclass[11pt]{amsart}

\title{Is every toric variety an M-variety?}

\author{Fr\'ed\'eric Bihan}
\address{Laboratoire de Math\'ematiques\\
         Universit\'e de Savoie\\
         73376 Le Bourget-du-Lac Cedex\\
         France}

\email{Frederic.Bihan@univ-savoie.fr}

\author{Matthias Franz}
\address{Fachbereich Mathematik \\ Universit\"at Konstanz \\ 
  78457 Konstanz \\ Germany}
\email{matthias.franz@ujf-grenoble.fr}

\author{Clint McCrory}
\address{Mathematics Department\\ University of Georgia
\\ Athens GA 30602 \\ USA}
\email{clint@math.uga.edu}

\author{Joost van Hamel}
\address{Departement Wiskunde\\ K.U. Leuven
\\ Celestijnenlaan 200B\\ B-3001 Leuven-Heverlee\\Belgium}
\email{vanhamel@member.ams.org}

\thanks{Research partially supported by NSF grant DMS-9810361 at  
MSRI}
\thanks{Bihan supported by F. Sottile via NSF CAREER grant DMS-0134860}
\thanks{Collaboration partially supported by
  Swiss National Science Foundation}

\numberwithin{equation}{section}

\theoremstyle{plain}
\newtheorem{thm}{Theorem}[section]

\newtheorem{prop}[thm]{Proposition}
\newtheorem{proposition}[thm]{Proposition}

\newtheorem{lemma}[thm]{Lemma}

\newtheorem{conjecture}[thm]{Conjecture}

\theoremstyle{definition}

\theoremstyle{remark}
\newtheorem{example}[thm]{Example}

\newtheorem{remark}[thm]{Remark}

\newcommand{\C}{\mathbb{C}}

\newcommand{\R}{\mathbb{R}}

\newcommand{\Z}{\mathbb{Z}}
\newcommand{\T}{\mathbb{T}}

\newcommand{\Or}{\mathcal{O}}
\newcommand{\clOr}{\overline{\mathcal{O}}}

\newcommand{\sI}{\mathcal{I}}

\DeclareMathOperator{\Gr}{Gr}
\DeclareMathOperator{\Ch}{Ch}

\newcommand{\Ztwo}{\Z/2}

\DeclareMathOperator{\Hom}{Hom}
\DeclareMathOperator{\Spec}{Spec}

\let\fcolon\colon

\def\cf{\textit{cf.}}

\begin{document}

\begin{abstract}
  A complex algebraic variety~$X$ defined over the real numbers
  is called an M-variety if the sum of its Betti numbers (for homology
  with closed supports and coefficients in~$\Ztwo$)
  coincides with the corresponding sum for the real part of~$X$.
  It has been known for a long time that any nonsingular complete
  toric variety is an M-variety. In this paper
  we consider whether this remains true for toric varieties that are  
  singular or not complete, and we give a positive answer
  when the dimension of~$X$ is less than or equal to~$3$.
\end{abstract}
\maketitle

\section{Introduction}
\noindent
Let $X$ be a topological space equipped with a continuous involution~$ 
\sigma$,
and let $X^{\sigma}$ denote
the fixed point set of~$\sigma$. For simplicity we assume that $X$ is a
finite-dimensional cell complex and $\sigma$ is a cellular involution.
The Smith--Thom inequality
asserts that the sum of the Betti numbers of~$X^{\sigma}$ with  
coefficients in~$\Ztwo$
does not exceed the the corresponding sum for~$X$,
\begin{equation}\label{E:Maximalinequality}
  \sum_k b_k(X^{\sigma}) \leq \sum_k b_k(X).
\end{equation}
Here and in the rest of the paper we consider ordinary homology groups,
or in the non-compact case homology groups with closed supports,  
also known as Borel--Moore homology.

A real algebraic variety is a complex variety~$X$
that is defined over the real numbers. Thus $X$ is
equipped with an antiholomorphic involution, complex conjugation. The
fixed point set of this involution is called the real part of~$X$ and  
will be
denoted by~$X(\R)$.  In contrast we shall often denote $X$ by~$X(\C)$.
The variety~$X$ is called an \emph{M-variety}
(\emph{maximal variety}) if equality occurs in~\eqref{E:Maximalinequality}.
In other words,
\begin{equation}\label{E:Maximalequality}
  \sum_k b_k(X(\R))= \sum_k b_k(X(\C)).
\end{equation}
M-varieties have attracted much attention in the study of
topological properties of real algebraic varieties.  One of the first
results in this domain is due to Harnack in~1876
(see \cite{Wilson} for a survey).
He proved an upper bound for the number of connected components of a
real algebraic curve in the projective plane; his theorem is a special
case of the Smith--Thom inequality.  He also showed that the bound is
sharp by constructing M-curves.  In the same vein, Itenberg and
Viro \cite{I-V} have recently constructed M-hypersurfaces of degree~$d$
in $n$-dimensional projective space for all positive integers $n$ and~$d$.

One of the most familiar M-varieties is projective space;
in this case the sums in~\eqref{E:Maximalequality} are equal to~$n+1$,
where $n$ is the dimension. Projective space is the simplest example
of a nonsingular complete toric variety.
Any toric variety is defined over the integers, hence over the reals.
There are many ways to see that every nonsingular complete toric variety
is an M-variety; see Section~\ref{sec:betti-complete}.
In fact, there is a degree-halving
isomorphism of algebras~$H^*(X(\C))\to H^*(X(\R))$
and the Betti numbers~$b_{2k}(X(\C))=b_{k}(X(\R))$
admit a simple combinatorial description in terms of fan data.


Another example of a toric M-variety is the complex algebraic torus itself,
whose real part is a real algebraic torus of the same rank.
Any toric variety is the disjoint union of torus orbits, hence of M-varieties.
In itself, the existence of
such a stratification is not enough for a variety to be an M-variety
(see for example Remark~\ref{rem:real-non-split}).
Nevertheless, our study of 
numerous examples has led us
to the following conjecture, which will be sharpened
in Section~\ref{sec:compare}.

\begin{conjecture}
  Every toric variety is an M-variety
  for homology with closed support.
\end{conjecture}

Note that for singular or non-complete toric varieties
there is no hope of getting a correspondence between individual Betti numbers:
In general the odd homology groups of~$X(\C)$ do not vanish,
and in particular the homology is not algebraic. In the nonsingular
non-complete case this already happens for the $1$-dimensional
torus; in the complete case such a phenomenon appears in
dimension~$2$ (see Proposition~\ref{surfaces} below, case~$(2)$).

Systematic studies of the stratification of a toric variety
by torus orbits and the associated
spectral sequence 
have been made by Fischli \cite{Fischli} and Jordan \cite{Jordan}.
By comparing this spectral sequence with another spectral sequence
for the homology of the real points, we arrive at our main result,
whose proof will appear in Section~\ref{sec:compare}.

\begin{thm}\label{main}
  Let $X$ be a (possibly singular or non-complete) toric variety
  of dimension~$\leq 3$.
  Then $X$ is an M-variety for homology with closed supports.
\end{thm}


In low dimensions straightforward spectral sequence calculations give
the individual 
Betti numbers of $X(\C)$ and~$X(\R)$;
see Proposition~\ref{surfaces}.


\section{Preliminaries
  }
\noindent
All vector spaces are over~$\Ztwo$ unless otherwise stated.
We write $H_*(X)$ for homology with \emph{closed supports}
(also known as \emph{Borel--Moore homology}) with coefficients in~$\Ztwo$.
Recall that for compact triangulable spaces, homology with closed
supports coincides with singular homology.

\medbreak

We now describe basic properties of toric varieties and review some
standard notation, referring to \cite{Fulton} and \cite{Oda} for details.

Any toric variety can be constructed in the following way:
Start with a lattice~$N$ of rank~$n$ and a rational fan~$\Delta$
in the real vector space~$N_{\R}=N \otimes_{\Z} \R$.
To each cone~$\sigma \in \Delta$ corresponds the
\emph{affine toric variety}~$X_{\sigma}=\Spec(\C[S_{\sigma}])$,
where $M=\Hom_{\Z}(N,\Z)$
is the lattice dual to~$N$, $M_{\R}=M \otimes_{\Z} \R$,
$\sigma^{\vee}=\{v \in M_{\R} \, | \, v(u) \geq 0 \; \forall \, u \in\sigma\}$
is the cone dual to~$\sigma$,
and $S_{\sigma}=\sigma^{\vee} \cap M$ is the corresponding semigroup.  
If $\tau$ is a face of~$\sigma$
, then
$X_{\tau}$ can be identified with a principal open subset of~$X_{\sigma}$.
The \emph{toric variety}~$X_{\Delta}$ is constructed by
gluing together the affine toric varieties~$X_{\sigma}$
along their common open subsets.
The (complex algebraic) torus associated with the lattice~$N$ is
$$
  \T_{N}:=\Spec(\C[M])=\Hom_{\Z}(M,\C^*)=N \otimes_{\Z} \C^*\simeq(\C^*)^n.
$$
The torus~$\T_{N}$ is open in~$X_{\Delta}$,
acts on each~$X_{\sigma}$, and this action extends to all of~$X_{\Delta}$
via the gluings.
The $\T_N$-orbits of~$X_{\Delta}$ are in one-to-one correspondence
with the cones in~$\Delta$ via the map~$\sigma \mapsto \Or_{\sigma}$,
where $\Or_{\sigma}$
is the $\T_N$-orbit of the distinguished point~$x_{\sigma}
\in X_{\sigma}=\Hom_{sg}(S_{\sigma},\C)$ defined by $x_{\sigma}(m)=1$ if
$-m \in S_{\sigma}$ and $x_{\sigma}(m)=0$ otherwise.
The Zariski closure~$\overline{\Or}_{\sigma}$
is the union of all orbits~$\Or_{\tau}$ such that
$\sigma$ is a face of~$\tau$.
For any cone~$\sigma$, define the lattices
$$
  N_{\sigma}:=(\sigma \cap N)+(-\sigma \cap N) \; ,
  \quad N 
  (\sigma):=N/N_{\sigma}.
$$
The lattice~$N_{\sigma}$ has rank~$\dim(\sigma)$ and
$N(\sigma)$ has rank~$n-\dim(\sigma)$.
The dual lattices~$M_{\sigma}=\Hom_{\Z}(N_{\sigma}, \Z)$
and~$M(\sigma)=\Hom_{\Z}(N(\sigma), \Z)$ are respectively
$$
  M_{\sigma} = M/(\sigma^{\bot} \cap M) \; ,
  \quad
  M(\sigma)=\sigma^{\bot} \cap M.
$$

The isotropy subgroup of the $\T_N$-action on~$\Or_{\sigma}$
is the isotropy group of the distinguished point~$x_{\sigma}$, which
consists of the~$t \in \T_N=\Hom_{\Z}(M,\C^*)$ such that $t(m)=1$
for any~$m \in \sigma^{\bot} \cap M$.
This gives the identification
$$
  \Or_{\sigma}=\T_{N(\sigma)}= \T_N / \T_{N_{\sigma}}.
$$


\section{Toric varieties as real varieties}
\noindent
This section is based on chapter 4 of \cite{Fulton}; see also \cite{Sottile}.
A toric variety~$X_\Delta$ is defined by polynomials with
integer coefficients. Thus we can consider $X_\Delta$ to be a real
algebraic variety, by which we mean a complex variety defined over the
real numbers. This is the standard real structure on a toric variety,
and it will be the only real structure we consider in this paper.

\begin{remark}\label{rem:real-non-split}
  An example of a toric variety with a non-standard real structure is
  the variety~$X={\C}P^1 \times{\C}P^1$ equipped with the
  antiholomorphic involution~$(z,w) \mapsto (\bar{w},\bar{z})$ (where
  the bar designates the usual complex conjugation on each factor).
  This type of real structure
(compatible with the action of a torus with a nonstandard involution)
has been studied in great detail by Delaunay
  \cite{Delaunay1}, \cite{Delaunay2}.  
  In this case the real part is homeomorphic to the $2$-sphere, and it is
  easy to check that $X$ is \emph{not} an M-variety.  
  On the other hand, taking
  the diagonal $D \subset {\C}P^1 \times{\C}P^1$, we have that both
  $D$ and its complement are M-varieties, so $X$ does admit a
  stratification where the open strata are M-varieties.
\end{remark}

As mentioned before, we will
denote the real part (i.e., the set of real points) of a real variety
$X$ by~$X({\R})$.
For clarity we shall often
denote the complex points by~$X(\C)$.
The real part $X_\Delta(\R)$ of a toric variety $X_\Delta$
is covered by the affine open subsets
$$
  X_{\sigma}(\R)=\Hom_{sg}(S_{\sigma},\R),
$$
where $\R$ is the multiplicative semigroup~$\R^* \cup \{0\}$, and
$$
  \T_{N}(\R)=\Spec({\R}[M])=\Hom_{\Z} (M, {\R}^*)
  = N \otimes_{\Z}\R^* \simeq(\R^*)^n,
$$
$$
  \Or_{\sigma}(\R)=\T_{N}(\R) \cdot x_{\sigma} = \T_{N(\sigma)}(\R)=
  \Hom_{\Z}(M(\sigma), \R^*).
$$
(Note that $x_\sigma\in X_{\sigma}(\R)$.)
The real part~$X_{\Delta}(\R)$ of a toric variety has
an orbit stratification similar
to that of the underlying complex toric variety.
$X_{\Delta}(\R)$ is obtained by gluing together
the~$X_{\sigma}(\R)$ for~$\sigma \in \Delta$, it is also the union of
the orbits~$\Or_{\sigma}(\R)\simeq {(\R^*)}^{n-\dim{\sigma}}$ under  
the action of~$\T_{N}(\R)$,
and the Zariski closure 
of~$\Or_{\sigma}(\R)$ is the union of all~$\Or_{\tau}(\R)$
such that
$\sigma$ is a face of~$\tau$.

As pointed out in \cite{Fulton}, this construction works for any
sub-semigroup of~$\C=\C^* \cup \{0\}$.
In particular, considering the semigroup~$\R_{+}=\R_{+}^* \cup \{0\}$
instead of~$\R=\R^* \cup \{0\}$, one
obtains the \emph{positive part} of a toric variety~$X_\Delta$. We will
denote the positive part of~$X_\Delta$ by~$X_\Delta^+$.
The positive part of~$X_\Delta$ is a semialgebraic subset of the real
part of~$X_\Delta$.

For any toric variety~$X_\Delta$,
we have $X_\Delta^+ \subset X_\Delta(\R) \subset X_\Delta(\C)$ due to the
semigroup inclusions~$\R_{+} \subset \R \subset \C$.
Moreover, the absolute value map~$z \rightarrow |z|$ gives rise to a
retraction~$\R_+ \subset \C \rightarrow \R_+$
which restricts to a retraction~$\R_+ \subset \R\rightarrow \R_+$.
The absolute value map can be extended in order to obtain the
following retractions.
$$
  \begin{array}{*5c}
    X_\Delta^{+} & \subset & X_\Delta(\C) & \rightarrow & X_\Delta^{+} \\
    \|            &         &    \cup        &            & \| \\
    X_\Delta^{+} & \subset & X_\Delta(\R) & \rightarrow & X_\Delta^{+}
  \end{array}
$$

For any lattice~$N$, define the \emph{compact torus}~$T_N$ by
$$
  T_N=\Hom_{\Z}(M,S^1) \subset \Hom_{\Z}(M,\C^*)=\T_N(\C),
$$
where $S^1$ is the unit circle in~$\C$.
Note that if $N$ has rank~$n$ then $T_N \simeq (S^1)^n$.
The set of $2$-torsion points of~$T_N$ will be
denoted by~$T_N[2]$. We have
$$
  T_N[2]=\Hom_{\Z}(M,S^0) \subset \Hom_{\Z}(M,\R^*)=\T_N(\R),
$$
where $S^0=\{\pm 1\}$ is the set of $2$-torsion points of~$S^1$.
If $N$ has rank~$n$, then $T_N[2] \simeq \{\pm 1\}^n$.
The isomorphism~$\C^* \simeq S^1 \times \R_+^*$ given by the map~$z
\mapsto
(z/ |z|,|z|)$ produces the identification
$$
  \T_N(\C) =\Hom_{\Z}(M,S^1) \times \Hom_{\Z}(M,\R_+^*)=T_N \times\T_N^+.
$$
Then, using the isomorphism~$\R_+^* \rightarrow \R$ given by the logarithm,
we obtain
$$
  \T_N^+=\Hom_{\Z}(M,\R_+^*)=\Hom_{\Z}(M,\R)=N_{\R},
$$
hence
$$
  \T_N(\C)=T_N \times N_{\R}.
$$
Similarly, we have that
$$
  \T_N(\R)=T_N[2] \times \T_N^+=T_N[2] \times N_{\R}.
$$
Applying this to the lattice~$N(\sigma)$ corresponding
to a cone~$\sigma \in\Delta$, we obtain
\begin{equation*}
  \Or_{\sigma}^+=N(\sigma)_{\R} \simeq \R^{n-\dim(\sigma)},\quad
  \Or_{\sigma}(\C)
  =T_{N(\sigma)}\times N(\sigma)_{\R}\;,
  \quad
  \Or_{\sigma}(\R)
  =T_{N(\sigma)}[2] \times N(\sigma)_{\R}\;.
\end{equation*}

From this discussion we get the following result, which is well-known
for the complex case (see \cite[\S 4.1, p. 79]{Fulton}).

\begin{prop}
  \label{mfdsing}
  The retraction~$r: X_{\Delta}(\C) \rightarrow  X_{\Delta}^{+}$ given
  by the absolute value map
  identifies~$X_{\Delta}^{+}$ with the quotient space of~$X_{\Delta}(\C)$
  by the action of the compact torus~$T_N$,
  and it also identifies $X_{\Delta}^{+}$ with the quotient space
  of~$X_{\Delta}(\R)$ by the action of~$T_N[2]$.
\end{prop}

The fibers of the quotient maps
$$
  X_{\Delta}(\C) \rightarrow X_{\Delta}^{+},
  \qquad
  X_{\Delta}(\R)\rightarrow X_{\Delta}^{+},
$$
over a point~$p \in X_{\Delta}^+$ are $T_{N(\sigma)}$ and $T_{N(\sigma)}[2]$,
respectively, where $\sigma \in \Delta$ is the unique cone
such that $p \in \Or_{\sigma}^+ \simeq \R^{n-\dim(\sigma)}$.
Since the exponential map gives the obvious identification
$$
  T_N[2]=({\textstyle\frac{1}{2}}N) /N \; \subset \; N_{\R} /N =T_N,
$$
where $\frac{1}{2}N=\{u \in N_{\R} \mid 2u \in N\}$,
Proposition \ref{mfdsing} can be restated as follows,
\cf~\cite[Theorem~11.5.4]{G-K-Z},~\cite[Proposition~4.1.1]{Delaunay2}.

\begin{prop}\label{mfdsing2}
  The toric variety~$X_{\Delta}(\C)$ is homeomorphic to the
  quotient space~$X_{\Delta}^+ \times T_N / \sim$,
  where the equivalence relation on~$X_{\Delta}^+ \times T_N$ is given
  by~$(p,t) \sim (p',t')$ if and only
  if $p=p'$ and $t-t' \in T_{N_{\sigma}}={(N_{\sigma})}_{\R} /{N_{\sigma}}$
  for the unique cone~$\sigma \in \Delta$ such that $p \in \Or_{\sigma}^+$.
  
  The real part~$X_{\Delta}(\R)$ is homeomorphic to the 
  quotient space~$X_{\Delta}^+ \times T_N[2] / \sim$,
  where the equivalence relation
  on~$X_{\Delta}^+ \times T_N[2]$ is given by~$(p,t) \sim (p',t')$ if and only
  if $p=p'$ and $t-t' \in T_{N_{\sigma}}[2]=\frac{1}{2}{N_{\sigma}} /{N_{\sigma}}$
  for the unique cone~$\sigma \in \Delta$ such that $p \in \Or_{\sigma}^+$.
\end{prop}


\section{Betti numbers of nonsingular complete toric varieties}
\label{sec:betti-complete}
\noindent
The fact that any nonsingular complete toric variety~$X=X_\Delta$
is an M-variety can be deduced from known results by the following arguments.

\begin{enumerate}
\item The Jurkiewicz--Danilov theorem~\cite[Proposition~10.4]{Danilov}
  implies that the
  cycle map from the Chow groups~$\Ch(X;R)$ to~$H_*(X;R)$
  is an isomorphism for arbitrary coefficients~$R$.
  Since the Chow groups are generated by closures of $\T_N$-orbits,
  which are conjugation-invariant subvarieties, it follows
  from standard results in equivariant cohomology
  (\cf~\cite[Remark~1.2.4\,(2)]{A-P}) that $X$ is maximal.
  Moreover, results of Krasnov~\cite{Krasnov} and
  Borel--Haefliger~\cite[\S 5.15]{Borel-Haefliger}
  imply that there exists a degree-halving isomorphism of algebras.
\item Using virtual Poincare polynomials (for homology
  with $\Ztwo$~coefficients), it is easy to show that the Betti numbers
  of~$X(\C)$ are the entries of the combinatorial $h$-vector
  of~$\Delta$, \cf~\cite[\S 4.5, \S 5.6]{Fulton} for the case
  of rational coefficients.
  One can imitate this proof for~$X(\R)$,
  using the virtual Poincar\'e polynomial 
  for real algebraic varieties
  defined by McCrory and Parusi\'nski~\cite{McCroryParusinski}.
  Here one starts with the virtual
  Poincar\'e polynomial~$\tilde P_{\R^*}(t) = t - 1$
  instead of~$\tilde P_{\C^*}(t) = t^2 - 1$.
  This implies the relations between the individual Betti numbers
  mentioned in the introduction:
  \begin{equation}
    \label{eq:relations-betti-h}
    b_{2k}(X_\Delta(\C))=b_k(X_\Delta(\R))=h_k(\Delta),
    \quad
    b_{2k+1}(X_\Delta(\C))=0.    
  \end{equation}
\item Nonsingular projective toric varieties are manifolds
  with a Hamiltonian torus action.
  Since all fixed points for the action on~$X(\C)$
  are contained in~$X(\R)$, a result of Duistermaat \cite[Theorem~3.1]{Du}
  implies that the Betti sum for~$X(\R)$
  is the number of these fixed points. Because the same is true for~$X(\C)$,
  this shows that every nonsingular projective toric variety is an M-variety.
Extending Duistermaat's methods, Biss, Guillemin and Holm also showed 
that there exists a degree halving isomorphism of 
algebras~$H^*(X(\C))\to H^*(X(\R))$ \cite[Corollary~5.8]{B-G-H}.

\item Relation~\eqref{eq:relations-betti-h} between the individual
  Betti numbers of the complex and real points of a nonsingular projective
  toric variety is also a special case of a result
  of Davis and Januszkiewicz~\cite[Theorem~3.1]{Davis-Janus}.
  Both Duistermaat's and Davis--Januszkiewicz's arguments are Morse-theoretic.
\item Another proof in the projective case is by ``shelling,''
  as in~\cite[\S 5.2]{Fulton}.
  (This is closely related to the 
  proof of Davis and Januszkiewicz.)
  If $X$ is nonsingular and projective, there is an ordering of the
  top-dimensional cones of~$\Delta$ that defines a filtration
  \begin{equation*}
    \emptyset \subset Z_m(\C)\subset\cdots\subset Z_1(\C) = X(\C)
  \end{equation*}
  by closed subvarieties~$Z_i(\C)$ with
  algebraic cells~$Y_i (\C)= Z_i(\C)\setminus Z_{i+1}(\C)\simeq\C^{k_i}$.
  The corresponding filtration~$Z_i(\R)$ of~$X(\R)$ has
  real algebraic cells~$Y_i(\R)\simeq \R^{k_i}$.
  Since the closure of the cell~$Y_i(\R)$ is a real algebraic variety,
  which is a cycle mod~$2$~\cite{Borel-Haefliger}, it follows
  that $b_{2k}(X(\C))= b_k(X(\R))$ and $b_{2k+1}(X(\C)) = 0$ for all~$k$.
\end{enumerate}




\section{The complex toric homology spectral sequence}
\label{sec:toric-hom-spsq}
\noindent
Following Totaro~\cite{Totaro}, Fischli~\cite{Fischli}, and
Jordan~\cite{Jordan},
we describe a spectral sequence associated to the filtration
of~$X=X_\Delta(\C)$ by the orbits of the torus action.
Let $\Delta^p$ be the set of cones of codimension~$p$
of the rational fan~$\Delta$ in~$N_\R \simeq \R^n$. Recall that for
each~$\sigma \in \Delta^p$ the orbit~$\Or_\sigma$ has dimension~$p$,
and the Zariski closure~$\clOr_\sigma$ is the union
of all orbits~$\Or_\tau$ such that $\sigma$ is a face of~$\tau$.
Taking
\begin{equation*}
  X_p := \bigcup_{\sigma \in \Delta^p} \clOr_\sigma
\end{equation*}
we get a filtration (in fact, a stratification)
\begin{equation}\label{filtration}
  \emptyset=X_{-1}\subset X_0\subset X_1\subset\dots\subset X_n = X_\Delta(\C)
\end{equation}
with each open stratum~$X_p^\circ := X_p \setminus X_{p-1}$
equal to the disjoint union of the $p$-dimensional orbits:
\begin{equation}\label{eq:disjoint-equidim-orbits}
  X_p^\circ = \bigcup_{\sigma \in \Delta^p} \Or_\sigma.
\end{equation}
Hence,
\begin{equation}\label{eq:orbit-homology-decomposition}
  H_i(X_p^\circ) = \bigoplus_{\sigma\in\Delta^p} H_i(\Or_\sigma).
\end{equation}

The filtration~\eqref{filtration} gives rise to a spectral sequence
\begin{equation}\label{eq:homology-spsq}
  E_{p,q}^1 
  = H_{p+q}(X_p^\circ) \Rightarrow H_{p+q}(X)
\end{equation}
converging to the homology (with closed supports) of~$X$
(\cf~\cite[p. 327]{MacLane}). The differentials
\begin{equation*}
  d_{p,q}^1 \fcolon
    H_{p+q}(X^\circ_p)
    \to
    H_{p+q-1}(X^\circ_{p-1})
\end{equation*}
at the $E^1$~level coincide with the connecting homomorphisms
in the long exact sequence
\begin{equation*}
  \dots \to
  H_{p+q}(X^\circ_{p-1}) \to
  H_{p+q}(X_p \setminus X_{p-2}) \to
  H_{p+q}(X^\circ_p) \to
  H_{p+q-1}(X^\circ_{p-1}) \to \dots
\end{equation*}
for the pair~%
$X^\circ_p \subset X_p \setminus X_{p-2}= X^\circ_p \cup X^\circ_{p-1}$.
It follows from \eqref{eq:orbit-homology-decomposition}
and from the functoriality of the connecting homomorphism
that we may describe each differential~$d_{p,q}^1$
in block matrix form~%
$(d_{q,\sigma, \tau})_{\sigma\in\Delta^p,\,\tau\in\Delta^{p-1}}$
with
\begin{equation*}
  d_{q,\sigma, \tau} \fcolon H_{p+q}(\Or_\sigma)\to H_{p+q-1}(\Or_\tau)
\end{equation*}
being the connecting homomorphism
for the pair~$\Or_\tau\subset\Or_\sigma\cup\Or_\tau$.
The latter is zero unless $\Or_\tau$ is in the boundary of~$\Or_\sigma$,
that is, unless $\sigma$ is a face of~$\tau$.

\begin{remark}\label{rem:leray-complex}
  The spectral sequence~\eqref{eq:homology-spsq} is
  isomorphic with the Leray spectral sequence of the
  retraction~$X_\Delta(\C)\to X^+_\Delta$. In particular,
  if $X_\Delta$ is projective, \eqref{eq:homology-spsq} is
  isomorphic with the Leray spectral
  sequence of the moment map.
\end{remark}

The differential~$d^1_{p,q}$ permits a simple description
in terms of the fan~$\Delta$, which we are going to derive now.
Recall that for each cone~$\sigma$ of codimension~$p$ we have
\begin{equation*}
  \Or_\sigma \simeq T_{N(\sigma)} \times \R^p.
\end{equation*}
Since
\begin{equation*}
  H_k(\R^{p}) =
  \begin{cases}
    \Ztwo & \text{if $k=p$,}\\
    0   & \text{otherwise,}
  \end{cases}
\end{equation*}
and since $\Ztwo$ has a unique generator, the K\"unneth formula gives
for any~$q$ a canonical isomorphism
\begin{equation}\label{eq:compare-homology}
   H_{p+q}(\Or_\sigma) = H_q(T_{N(\sigma)})
\end{equation}
in homology with closed supports.
Hence each component~$d_{q,\sigma, \tau}$ of the differential~$d^1$
is canonically identified with a map
\begin{equation*}
d'_{q,\sigma, \tau} \fcolon H_{q}(T_{N(\sigma)})
\to H_{q}(T_{N(\tau)}).
\end{equation*}
(If the characteristic of the coefficients was different from~$2$,
we would have to choose orientations, and signs would appear
as in simplicial homology.)

If $\sigma$ is a face of~$\tau$, then the
lattice~$N_\sigma \subset N$ of the isotropy group of~$\T_N$ acting
on~$\Or_\sigma$ is contained in the the corresponding
lattice~$N_\tau$.
Hence we get a natural
surjection~$N(\sigma) = N/N_\sigma\to N/N_\tau = N(\tau)$,
which gives natural split surjections~$\T_{N(\sigma)} \to \T_{N(\tau)}$
and $T_{N(\sigma)} \to T_{N(\tau)}$,
both of which we will denote by~$\pi_{\sigma, \tau}$.


\begin{proposition}[{\cite[Thm~2.1]{Fischli}, \cite[\S 2.3]{Jordan}}]
\label{prop:Jordan}
  If $\sigma$ is a facet of~$\tau$, then the homomorphism
  \begin{equation*}
    d'_{q,\sigma, \tau} \fcolon H_{q}(T_{N(\sigma)})
    \to H_{q}(T_{N(\tau)}).
  \end{equation*}
  coincides with the homomorphism induced by the split surjection
  \begin{equation*}
    \pi_{\sigma, \tau} \fcolon T_{N(\sigma)} \to T_{N(\tau)}.
  \end{equation*}
\end{proposition}

\begin{proof}
  See for example \cite[\S 2.3]{Jordan}, where it is derived from
  Proposition~\ref{mfdsing}.
  Alternatively, it follows easily from the fact that
  the pair~$\Or_\tau\subset\Or_\sigma\cup\Or_\tau$ is isomorphic with
  the pair~$(\C^*)^{p-1}\times\{0\}\subset(\C^*)^{p-1}\times\C$
  for~$\sigma\in\Delta^p$. To see this, consider the star of~$\sigma$,
  that is, the set of cones~$\tau\in\Delta$ having $\sigma$ as face.
  Taking the image of each such~$\tau$ in~$N(\sigma)$ gives a new fan,
  whose associated toric variety is~$\clOr_\sigma$,
  \cf.~\cite[pp.~52--54]{Fulton}.
  Hence, the pair~$(\Or_\tau,\Or_\sigma)$ is isomorphic with the pair
  given by ray and the origin in~$N(\sigma)$.
\end{proof}

Totaro~\cite{Totaro} noticed that for rational coefficients
this spectral sequence degenerates at the $E^2$~level
for any toric variety for
any toric variety (\cf~\cite[Proposition~2.4.5]{Jordan}).
In~\cite{Franz2}, the second author proved that
over an arbitrary coefficient ring~$R$ this happens
for toric varieties which are $R$-homology manifolds,
in particular in the nonsingular case.
Moreover, he conjectured that degeneration occurs for
arbitrary toric varieties, as with rational coefficients.


\section{Homology of compact tori and their 2-torsion points}
\label{sec:hom-comp-tori}
\noindent
In order to construct a spectral sequence for~$X(\R)$ which compares well
with the spectral sequence introduced for~$X(\C)$,
we need to relate the homology of an $n$-dimensional compact torus~$T$
and its $2$-torsion points~$T[2]$.

The homology of any compact topological group~$G$ is a graded algebra
by the \emph{Pontryagin product}~$H_*(G) \otimes H_*(G) \to H_*(G)$,
which is constructed in the obvious way from the
Eilenberg--Zilber map~$H_*(G) \otimes H_*(G)\to H_*(G \times G)$
and the multiplication~$G \times G\to G$.
(The restriction to compact groups is caused by our choice of
homology with closed supports.)
The unit in~$H_*(G)$ is the homology class~$[1]$
of the identity element~$1\in G$.
Since $T$ and $T[2]$ are commutative, the Pontryagin product
is (anti)commutative in these cases.
Note that $H_*(T[2])=H_0(T[2])$ is nothing but the group algebra of the finite
group~$T[2]$ with coefficients in~$\Ztwo$.

$H_*(S^1)$ is an exterior algebra on
the fundamental class of~$S^1$. Similarly, the homology
of~$S^1[2]=S^0=\{1,g\}$ is an exterior algebra on the fundamental
class~$[1]+[g]$ because
$$
  \bigl([1]+[g]\bigr)^2=[1]^2+[g]^2=2\,[1]=0.
$$
Hence, $H_*(S^1)$ and $H_0(S^0)$ are isomorphic as ungraded algebras.
Note that in both cases the generators are unique, so that the isomorphism
is canonical.

Choosing an decomposition of~$T$ into circles
\begin{align}
   \label{eq:T-product}
   T &\simeq (S^1)^n,
\intertext{the K\"unneth formula gives isomorphisms}
   \label{eq:T-Kuenneth}
   H_*(T) &\simeq H_*(S^1)^{\otimes n}, \\
   \label{eq:T2-Kuenneth}
   H_0(T[2]) & \simeq H_0(S^0)^{\otimes n}.
\end{align}
Since \eqref{eq:T-product} is a decomposition
of topological groups, \eqref{eq:T-Kuenneth}~and~\eqref{eq:T2-Kuenneth}
are isomorphisms of 
algebras.
This gives an isomorphism
\begin{equation}
  \label{eq:iso-torus}
  H_*(T)\simeq H_0(T[2]) \quad\hbox{as ungraded algebras.}
\end{equation}

The following example shows that this isomorphism does depend
on the decomposition of~$T$ as a product of copies of~$S^1$ and that
it is not functorial, except for homomorphisms that are compatible  
with the chosen product decompositions.

\begin{example}\label{ex:class-subtorus}
  Let $T=S^1 \times S^1$, and let $D \subset T$ be the diagonal torus. Then
  \begin{align*}
    [D] & = [S^1\times1] + [1\times S^1] \in
      H_*(T),
  \intertext{whereas}
    [D[2]] & =  [S^0\times 1] + [1\times S^0] + [S^0\times S^0] \in
      H_0(T[2]).
  \end{align*}
  By functoriality, the same relations hold for any two circles
  representing different elements in~$H_*(T)$ and a circle
  representing their sum.
\end{example}

\medskip

In order to get a better comparison of $H_*(T)$ and $H_0(T[2])$
we will construct a natural filtration on~$H_0(T[2])$.
Let $\sI=\sI(T[2])$ be the kernel of the
augmentation~$H_0(T[2]) \to \Ztwo$
induced by the group homomorphism from $T$ to a point.
It is an ideal in~$H_0(T[2])$, and
there is a canonical direct sum decomposition of vector spaces
\begin{equation}\label{eq:decomposition}
  H_0(T[2]) \simeq \langle [1] \rangle \oplus \sI,
\end{equation}
functorial with respect to homomorphisms to tori.

\begin{lemma}
  For~$k\ge1$, the $k$-th power~$\sI^k$ of~$\sI$ is additively
  generated by the fundamental classes of rank~$k$~subtori of~$T[2]$.
  In particular, $\sI^k=0$ if $k$ is greater than the rank~$n$ of~$T[2]$.
\end{lemma}

\begin{proof}
  An element~$a\in H_0(T[2])$ is a formal linear combination of points
  of~$T[2]$. We have $a\in\sI$ if and only if there is an even number
  of such points.
  Denoting by~$A$ the set of these points, we can rewrite $a$ as
  $$
    a=\sum_{g\in A}[g]=\sum_{1\ne g\in A}\bigl([1]+[g]\bigr),
  $$
  which is a sum of fundamental classes of rank~$1$~subtori.

  Suppose the claim is true for~$k$. Then $\sI^{k+1}$ is additively generated
  by products~$[H]*([1]+[g])$ with~$g\in T$ and $H$ a rank~$k$~subtorus.
  If $g$ is contained in~$H$, this equals $[H]+[H]=0$. Otherwise we get
  the fundamental class of a subtorus of rank~$k+1$.
\end{proof}

We now look at graded quotient associated to the filtration
$$
  H_0(T[2])=\sI^0\supset\sI^1\supset\dots\supset\sI^n\supset \sI^{n+1} = 0.
$$

Writing $V = \Hom(S^1, T) \otimes \Ztwo$,
the assignment~$\chi \mapsto \chi_*([S^1])$ gives a canonical isomorphism
of vector spaces~$V = H_1(T)$, which extends to an isomorphism
of graded algebras
\begin{equation}\label{iso-V-H}
  \bigwedge\nolimits^{\!*} V = H_*(T),
\end{equation}
natural with respect to homomorphisms of tori.

\begin{proposition}\label{prop:graded-comparison}
   There is a natural isomorphism of graded algebras
   \begin{equation*}
     \Gr_\sI^* H_0(T[2]) = H_*(T).
   \end{equation*}
\end{proposition}

\begin{proof}
  The assignment~$\chi \mapsto \chi_*([S^0])$
  gives a map (not a homomorphism)~$V \to \sI \subset H_0(T[2])$,
  which induces a map
  \begin{equation*}
    V \simeq \sI/\sI^2.
  \end{equation*}
  We claim that the latter is a homomorphism, i.e., compatible with sums:
  Take two elements from~$V$. We may assume that they are distinct and
  non-zero. Then Example~\ref{ex:class-subtorus} shows that the image
  of their sum and the sum of their images in~$\sI$
  differ by an element from~$\sI^2$.

  The square of any element~$a\in\sI$ is zero because $a$ is a sum
  of terms of the form~$[1]+[g]$.
  Therefore, we get a homomorphism of graded algebras
  \begin{equation}\label{iso-V-Gr}
    \bigwedge\nolimits^{\!*} V \to \Gr_\sI^* H_0(T[2]).
  \end{equation}
  Since $\sI^k$ is additively generated by the fundamental classes
  of rank~$k$ subgroups, the map~\eqref{iso-V-Gr} is surjective,
  hence an isomorphism.

  The desired map is the composition of \eqref{iso-V-H}~and~\eqref{iso-V-Gr}.
  It is clear from the definitions that both isomorphisms do not
  depend on any choice 
  and that they are functorial in~$T$.
\end{proof}

\begin{remark}\label{rem:naturality-iso-T2}
  Together with the isomorphism~\eqref{eq:iso-torus},
  Proposition~\ref{prop:graded-comparison} implies
  that $H_0(T[2])$ and $\Gr_\sI^* H_0(T[2])$ are isomorphic
  as algebras. Again the isomorphism is not natural in general,
  but it is so with respect to homomorphisms that are compatible  
  with the chosen product decompositions.
  In particular, one can choose isomorphisms compatible with
  a single injection or a single split projection.
\end{remark}


\section{The real toric homology spectral sequence}
\label{sec:real-spectral-sequence}
\noindent
As for~$X(\C)$, one can consider for the real toric variety~$X(\R)$
the filtration by $\T_N(\R)$-orbit dimension.
Because $\T_{N(\sigma)}=T_{N(\sigma)}[2]\times\R^p\simeq(\R^*)^p$
is a disjoint union of $p$-cells
for~$\sigma\in\Delta^p$, the resulting filtration actually is cellular.
This leads to a cellular chain complex~$C_*(\Delta)$, and the associated
spectral sequence, which is again isomorphic with the Leray spectral
sequence of~$X_\Delta(\R)\to X^+_\Delta$, degenerates at the $E^1$~level.

As in the complex case, we can identify
$$
  C_p(\Delta)=\bigoplus_{\sigma\in\Delta^p} H_0(T_{N(\sigma)}[2]),
$$
and we have the following result.

\begin{lemma}\label{lem:real-toric-hom}
  Under this identification,
  the component~$H_0(T_{N(\sigma)}[2])\to H_0(T_{N(\tau)}[2])$
  of the differential 
  is induced by the surjection~$T_{N(\sigma)}[2]\to T_{N(\tau)}[2]$
  for $\sigma$ a facet of~$\tau$, and zero otherwise.
\end{lemma}

Using the decreasing filtration on $2$-tori introduced in the previous
section, we define an increasing filtration on~$C_*(\Delta)$:
$$
  F_p(C_*(\Delta))=\bigoplus_{\sigma\in\Delta}
                     \sI^{-p}(T_{N(\sigma)}[2])
  \quad\hbox{for~$p\le0$.}
$$
Then $F_{-n-1}(C_*(\Delta))=0$ and $F_0(C_*(\Delta))=C_*(\Delta)$.
We write $G^k_{pq}$ for the terms of the resulting spectral sequence,
which is concentrated in the second quadrant. For~$n=3$, the possibly
non-zero terms are located as follows:
$$
  \begin{array}{*4{|c}|l}
    \cline{1-4}
    * &   &   &   & 6=q \\
    \cline{1-4}
      & * &   &   &   \\
    \cline{1-4}
      & * & * &   &   \\
    \cline{1-4}
      &   & * & * &   \\
    \cline{1-4}
      &   & * & * &   \\
    \cline{1-4}
      &   &   & * &   \\
    \cline{1-4}
      &   &   & * & 0 \\
    \cline{1-4}
    \multicolumn1r{\hbox to0pt{\hss$p=-3$}} &
      \multicolumn2c{} & \multicolumn1c 0
  \end{array}
$$

As a direct consequence of Proposition~\ref{prop:graded-comparison},
we obtain the following result.

\begin{proposition}\label{prop:compare-spectral-sequences}
  The term~$G^0_{pq}$ of the spectral sequence for $X(\R)$ is isomorphic
  to the term~$E^1_{p+q,-p}$ of the spectral sequence for~$X(\C)$,
  and this is compatible with the differentials.
  Hence also $G^1_{pq}=E^2_{p+q,-p}$. In particular, the $G^1$~and~$E^2$~levels
  have the same total dimension.
\end{proposition}

\begin{proof}
By construction, $G^0_{pq} = \bigoplus_{\sigma \in \Delta^{p+q}}
\Gr_\sI^{-p} H_0(T_{N(\sigma)}[2])$.
It now follows from 
Proposition~\ref{prop:graded-comparison} and the discussion in
Section~\ref{sec:toric-hom-spsq} that $G^0_{pq}= E^1_{p+q,-p}$ for all
$p,\,q$. The fact that for each $p,\,q$ differential
$d_{pq}^0 \fcolon G^0_{pq} \to G^0_{p,q-1}$ corresponds to the
differential
$d_{pq}^1 \fcolon E^1_{p+q, -p} \to E^1_{p+q-1,-p}$ is easy to check
from
Lemma~\ref{lem:real-toric-hom} and Proposition~\ref{prop:Jordan}.
\end{proof}

\begin{remark}\label{rem:spseq-diff}
The isomorphism $G^1_{pq} \simeq E^2_{p+q, -p}$ is not compatible with the
differentials. Indeed, the differential with source $G^1_{pq}$ has
target $G^1_{p-1,q} \simeq E^2_{p+q-1, -p+1}$, 
whereas the differential with source
$E^2_{p+q,-p}$ has target $E^2_{p+q-2,-p+1}$.
\end{remark}

\begin{remark}\label{real-split}
  The decomposition~\eqref{eq:decomposition} induces a
  direct sum decomposition of~$C_*(\Delta)$. One summand is
  the complex~$G^0_{0,*}=E^1_{*,0}$, the other one the direct sum of
  the augmentation ideals~$\sI(T_{N(\sigma)})$.
  Therefore, all differentials starting at the rightmost column~$G^k_{0,*}$
  vanish for~$k\ge1$.
\end{remark}


\section{Comparing the spectral sequences}\label{sec:compare}
\noindent
We can now sharpen the conjecture stated in the introduction.

\begin{conjecture}\label{conj:degeneration}
  The spectral sequence~$G^k_{p,q}$ for~$H_*(X(\R))$
  degenerates at the $G^1$~level.
\end{conjecture}

Notice that this implies that the spectral sequence~$E^k_{p,q}$
converging to~$H_*(X(\C))$ degenerates at the $E^2$~level: the
comparison between $G^1$ and $E^2$
(Proposition~\ref{prop:compare-spectral-sequences}) implies that 
otherwise the
Betti sum for~$X(\R)$ would exceed that for~$X(\C)$, in
contradiction to the Smith--Thom
inequality~\eqref{E:Maximalinequality}.
In particular, the conjecture implies that $X$ is an M-variety.

Numerous examples indicate that degeneration is at least
not a rare phenomenon. We have checked many examples 
using the Maple package Torhom by the second author~\cite{Franz}.
So far this has
not led to an example of a toric variety for which
the spectral sequence~$G^k_{pq}$ does not degenerate.
Here is one example:

\begin{example}
  Let $\Delta$ be the normal fan (\cf~\cite[Section~5.6]{Fulton})
  of the five-dimensional cyclic polytope
  with vertices~$(k^1,\ldots,k^5)$ for~$k=0$,~$1$,~\ldots,~$6$.
  The dimensions of the $E^2$~terms and the $G^1$~terms are:
  $$
    \begin{array}[b]{r|*6{r|}}
      \cline{2-7}
      q=5 & \phantom{00} & \phantom{00} & \phantom{00} &
          \phantom{00} & \phantom{00} & \phantom{0}1 \\
      \cline{2-7}
        &    &    &    &    & 11 &  4 \\
      \cline{2-7}
        &    &    &    & 13 & 27 &  6 \\
      \cline{2-7}
        &    &    &  6 & 17 & 21 &  4 \\
      \cline{2-7}
        &    &  1 &  1 &  4 &  5 &  1 \\
      \cline{2-7}
      0 &  1 &  0 &  0 &  0 &  0 &  0 \\
      \cline{2-7}
      \multicolumn1c{} & \multicolumn1r 0 & \multicolumn4c{} &
        \multicolumn1l{\hbox to0pt{$\phantom{0}5=p$\hss}} \\
      \multicolumn7c{\mskip60mu E^2_{pq}}
    \end{array}
    \qquad
    \qquad
    \begin{array}[b]{r|*6{r|}l}
      \cline{2-7}
      & \phantom{0}1 & \phantom{00} & \phantom{00} &
          \phantom{00} & \phantom{00} & \phantom{00} & q=10 \\
      \cline{2-7}
      &    &  4 &    &    &    &    &    \\\cline{2-7}
      &    & 11 &  6 &    &    &    &    \\\cline{2-7}
      &    &    & 27 &  4 &    &    &    \\\cline{2-7}
      &    &    & 13 & 21 &  1 &    &    \\\cline{2-7}
      &    &    &    & 17 &  5 &  0 &    \\\cline{2-7}
      &    &    &    &  6 &  4 &  0 &    \\\cline{2-7}
      &    &    &    &    &  1 &  0 &    \\\cline{2-7}
      &    &    &    &    &  1 &  0 &    \\\cline{2-7}
      &    &    &    &    &    &  0 &    \\\cline{2-7}
      &    &    &    &    &    &  1 &   0\\\cline{2-7}
        \multicolumn1c{\hbox to0pt{$p=-5$\hss}} & 
        \multicolumn5c{} & \multicolumn1c 0 \\
      \multicolumn7c{\mskip30mu G^1_{pq}}
    \end{array}
  $$
  Clearly there is plenty of room for higher differentials in both cases.
  Yet both spectral sequences degenerate, as Torhom's direct computation
  of~$H_*(X(\R))$ shows.
\end{example}

We now prove our main result 
by verifying Conjecture~\ref{conj:degeneration}
for varieties up to dimension~$3$.

\begin{proof}[Proof of Theorem~\ref{main}]
  There can be no higher differentials starting
  at~$G^1_{0,*}$~ by Remark~\ref{real-split}. Hence the only
  higher differential can be~$d^1\fcolon G^1_{-1,4}\to G^1_{-2,4}$
  for $3$-dimensional varieties. 
  Clearly, the source of this differential
  is non-trivial only if the previous differential
  $d^0\fcolon G^0_{-1,4} \to G^0_{-1,3}$ is not injective.
  By Proposition~\ref{prop:compare-spectral-sequences} this is
  the case only if the differential $d^1 \fcolon E^1_{3,1} \to E^1_{2,1}$ 
  is not injective, in other words, if the map
  $$
    H_1(T_N)\to\bigoplus_{\tau\in\Delta^2} H_1(T_{N(\tau)})
  $$
  is not injective. Since the kernel of each map~$H_1(T_N)\to H_1(T_{N(\tau)})$
  is generated by the image of the minimal representative of~$\tau$
  in~$V=N/2N$, this implies that all $\tau\in\Delta^2$ have the same
  image~$v$ in~$V$. As a consequence, all 2-tori
  $$
    T_{N(\sigma)}[2]=({\textstyle\frac12}N/v)\!\bigm/\!(N/v)=:\tilde T
  $$
  of rank~$2$ are actually identical.
  By Remark~\ref{rem:naturality-iso-T2}, one can therefore choose
  isomorphisms
  $H_0(T_N[2])\simeq\Gr_\sI^* H_0(T_N[2])$ and
  $H_0(\tilde T)\simeq\Gr_\sI^* H_0(\tilde T)$
  which intertwine the differential
  $$
    d^0\fcolon
      \bigoplus_{p+q=3}G^0_{p,q}=\Gr_\sI^* H_0(T_N[2])\to
      \bigoplus_{\sigma\in\Delta^2}\Gr_\sI^* H_0(T_{N(\sigma[2])})
      =\bigoplus_{p+q=2}G^0_{p,q}
  $$
  and the underlying differential
  $$
    d\fcolon C_3(\Delta)=H_0(T_N[2])\to
      \bigoplus_{\sigma\in\Delta^2}H_0(T_{N(\sigma[2])})=C_2(\Delta).
  $$
  Hence the kernels of these differentials have the same dimension,
  and all higher differentials~$d^k_{p,q}$ must vanish
  for~$p+q=3$ and~$k\ge1$.
\end{proof}

Here is another case where the spectral sequences degenerates.

\begin{proposition}
  Conjecture~\ref{conj:degeneration} is true
  for complete toric varieties with isolated singularities.
\end{proposition}

\begin{proof}
  By~\cite[\S1.2]{Brion}, the term~$E^2_{pq}$ is concentrated on the two lines
  $p=q$ and $p=q+1$ because only full-dimensional cones may not be regular.
  (In~\cite{Brion} rational coefficients are used.
  The same arguments work for~$\Ztwo$ if one considers regular
  (or $\Ztwo$-regular) cones instead of simplicial ones.)
  Therefore, non-trivial higher differentials are impossible
  for~$G^1_{pq}$.
\end{proof}


\section{Explicit calculations in dimension~2}
\noindent
Jordan \cite[Theorem~3.4.2]{Jordan} computed
the integral homology with closed supports of an arbitrary
2-dimensional toric variety.
Together with our result, this leads to the following classification
in the case of complete toric surfaces.

\begin{proposition}\label{surfaces} Let $X$ be the complete toric surface
   associated to a fan~$\Delta$ with respect to a lattice~$N$ of rank~$2$.
   Let $s$ be the number of $1$-dimensional cones of~$\Delta$.

   \begin{itemize}

   \item[(1)] If at least two primitive generators of one-dimensional cones
     of~$\Delta$ have different images in the quotient lattice~$N/ 2N$, then
     $$
       b_0(X(\R))=1 \, , \; b_1(X(\R))=s-2 \, , \; b_2(X(\R))=1,
     $$
     $$
       b_0(X(\C))=1 \, , \; b_1(X(\C))=0 \, , \; b_2(X(\C))=s-2 \, ,
       \;b_3(X(\C))=0\, ,  \; b_4(X(\C))=1.
     $$
   \item[(2)] If the primitive generators of one-dimensional
     cones of~$\Delta$ all have the same image in~$N/ 2N$, then
     $$
       b_0(X(\R))=1 \, , \; b_1(X(\R))=s-1 \, , \; b_2(X(\R))=2,
     $$
     $$
       b_0(X(\C))=1 \, , \; b_1(X(\C))=0 \, , \; b_2(X(\C))=s-1 \, ,
       \; b_3(X(\C))=1 \, , \; b_4(X(\C))=1.
     $$
\end{itemize}
\end{proposition}

\begin{proof}
  Recall that every complete two-dimensional fan~$\Delta$
  is the normal fan of some polytope~$P$.
  The dimensions of the $E^1$~level of the spectral sequence
  for~$X(\C)$ are as follows,
  where the $q=0$~row is the cellular chain complex of such a~$P$:
  $$
    \begin{array}{r|*3{c|}}
      \cline{2-4}
      q=2 & \makebox[3em]{}  & \makebox[3em]{}  & \makebox[3em]{$1$} \\
      \cline{2-4}
        &   & s & 2 \\
      \cline{2-4}
      0 & s & s & 1 \\
      \cline{2-4}
      \multicolumn1c{} & \multicolumn1c 0 & \multicolumn1l{} &
        \multicolumn1c{\!\!\!\hbox to0pt{$2=p$\hss}}
    \end{array}
  $$

  The $q=0$~row of~$E^2$ is the homology of~$P$.
  As in the proof of Theorem~\ref{main}, we know that the
  differential~$E^1_{2,1}\to E^1_{1,1}$
  is either injective or has a one-dimensional kernel. The latter occurs
  if and only if the minimal generators of the one-dimensional
  cones in~$\Delta$ have the same image in~$N/2N$.
  Hence, the $E^2$ term is one of the following:
  $$
    (1)
    \begin{array}{r|*3{c|}}
      \cline{2-4}
      q=2 & \makebox[3em]{} & \makebox[3em]{}  &  \makebox[3em]{$1$} \\
      \cline{2-4}
        &   & s-2 & 0 \\
      \cline{2-4}
      0 & 1 & 0 & 0 \\
      \cline{2-4}
      \multicolumn1c{} & \multicolumn1c 0 & \multicolumn1l{} &
        \multicolumn1c{\!\!\!\hbox to0pt{$2=p$\hss}}
    \end{array}
     \qquad\qquad
    (2)
    \begin{array}{r|*3{c|}}
      \cline{2-4}
      q=2 & \makebox[3em]{}  & \makebox[3em]{}  & \makebox[3em]{$1$} \\
      \cline{2-4}
        &   & {s-1} & 1 \\
      \cline{2-4}
      0 & 1 & 0 & 0 \\
      \cline{2-4}
      \multicolumn1c{} & \multicolumn1c 0 & \multicolumn1l{} &
        \multicolumn1c{\!\!\!\hbox to0pt{$2=p$\hss}}
    \end{array}
  $$
  From this we can read off the Betti numbers of~$X(\C)$.
  To get the Betti numbers of~$X(\R)$, we simply
  have to sum up the dimensions in each column.
\end{proof}

If $X$ is a nonsingular complete toric surface, then the primitive
generators of two consecutive (with respect to either of the two cyclic orders)
one-dimensional cones have different images in $N/2N$, so that
case~(1) applies.
It is instructive to use Proposition~\ref{mfdsing2} to find
generators for the homology groups of $X(\R)$.


\bibliographystyle{amsalpha}

\begin{thebibliography}{BGH}

\bibitem[AP]{A-P}
  C.~Allday and V.~Puppe,
  \textit{Cohomological methods in transformation groups},
  Cambridge University Press, Cambridge 1993.



\bibitem[BGH]{B-G-H}
  D. Biss, V. Guillemin, T. Holm,
  \emph{The mod 2 cohomology of fixed point sets of
    anti-symplectic involutions},
  Adv. Math. {\bf 185} (2004), 370--399.

\bibitem[BH]{Borel-Haefliger}
  A. Borel and A. Haefliger,
  \textit{La classe d'homologie fondamentale d'un espace analytique},
  Bull. Soc. Math. France \textbf{89} (1961), 461--513.



\bibitem[Br]{Brion}
  M. Brion,
  \textit{The structure of the polytope algebra},
  T\^ohoku Math. J. \textbf{49} (1997), 1--32.


\bibitem[D]{Danilov}
  V. I. Danilov,
  \textit{The geometry of toric varieties},
  Russian Math. Surveys \textbf{33} (2) (1978), 97--154.

\bibitem[DJ]{Davis-Janus} M. Davis and T. Januszkiewicz, \textit{Convex
polytopes, Coxeter orbifolds and torus actions}, Duke Math. J. \textbf 
{62} (1991),
417--451.

\bibitem[De1]{Delaunay2}
  C. Delaunay,
  \textit{Real structures on smooth compact toric surfaces},
  Topics in Algebraic Geometry and Geometric Modeling, R. Goldman and R.
  Krasuaskas eds., Contemp. Math. \textbf{334}, Amer. Math. Soc.,
  Providence, RI, 2003, 267--290.

\bibitem[De2]{Delaunay1} C. Delaunay,
  \textit{Real structures on compact toric varieties},
  Ph.\,D. thesis, Universit\'e de Strasbourg, 2004.

\bibitem[Du]{Du}
  H. Duistermaat,
  \emph{Convexity and tightness for restrictions of Hamiltonian  
    functions to fixed point sets of an antisymplectic involution},
  Trans. Amer. Math. Soc. {\bf 275} (1) (1983), 417--429.

\bibitem[Fi]{Fischli} S. Fischli, \textit{On toric varieties},
  Ph.\,D. thesis, Universit\"at Bern, 1992. \\
  Available at \texttt{http://www.hta-be.bfh.ch/\~{ }fischli/}

\bibitem[Fr1]{Franz}
  M. Franz,
  Maple package \texttt{Torhom}, version 1.3.0, September 13, 2004.\\
  Available at
  \texttt{http://www-fourier.ujf-grenoble.fr/\~{ }franz/maple/torhom.html}

\bibitem[Fr2]{Franz2} M.~Franz,
  \textit{The integral cohomology of smooth toric manifolds},
  to appear in Proc. Steklov Inst. Math.

\bibitem[F]{Fulton}
  W. Fulton,
  \textit{Introduction to Toric Varieties},
  Princeton University Press, Princeton 1993.

\bibitem[GKZ]{G-K-Z}
  I. Gelfand, V. Kapranov, and A. Zelevinsky,
  \textit{Discriminants, Resultants, and Multidimensional Determinants},
  Birkh\"auser, Boston, 1994.

\bibitem[IV]{I-V}
  I. Itenberg and O. Viro,
  \textit{Maximal real algebraic hypersurfaces of projective space},
  to appear.

\bibitem[J]{Jordan}
  A. Jordan,
  \textit{Homology and cohomology of toric varieties},
  Ph.\,D. thesis, Universit\"at Konstanz,
  Konstanzer Schriften in Mathematik und Informatik
  \textbf{57}, 1998. \\
  Available at
  \texttt{http://www.inf.uni-konstanz.de/Schriften/preprints-1998.html\#057}

\bibitem[Kr]{Krasnov}
V.A. Krasnov, \textit{Real algebraically maximal varieties},
Mat. Zametki, \textbf{73} (2003), 853--860,
  English translation in:  Math. Notes  \textbf{73}  (2003), 806--812.

\bibitem[M]{MacLane}
  S. Mac Lane,
  \textit{Homology},
  Springer Verlag, New York, 1967.

\bibitem[MP]{McCroryParusinski} C. McCrory and A. Parusi\'nski,
  \textit{Virtual Betti numbers of real algebraic varieties},
  Comptes Rendus Acad. Sci. Paris, Ser.~I \textbf{336} (2003), 763--768.

\bibitem[O]{Oda}
  T. Oda,
  \textit{Convex Bodies and Algebraic Geometry},
  Springer-Verlag, New York, 1988.

\bibitem[S]{Sottile}
  F. Sottile,
  \textit{Toric ideals, real toric varieties, and the moment map},
  Topics in Algebraic Geometry and Geometric Modeling, R. Goldman and R.
  Krasuaskas eds., Contemp. Math. \textbf{334}, Amer. Math. Soc.,
  Providence, RI, 2003, p. 225--240.

\bibitem[T]{Totaro}
  B. Totaro,
  \textit{Chow groups, Chow cohomology, and linear varieties},
  to appear in J.\ Alg.\ Geom. \\
  Available at \texttt{http://www.dpmms.cam.ac.uk/\~{ }bt219/papers.html}

\bibitem[W]{Wilson} G. Wilson, \textit{Hilbert's sixteenth problem},
Topology
\textbf{17} (1978), 53--73.

\end{thebibliography}

\end{document}